\documentclass{amsart}
\usepackage{amssymb,amsthm,multicol,hyperref}
\PassOptionsToPackage{unicode}{hyperref}
\PassOptionsToPackage{naturalnames}{hyperref}
\newtheorem{thm}{Theorem}
\newtheorem{lem}{Lemma}
\newtheorem*{thm*}{Theorem}
\newtheorem{cor}{Corollary}
\newtheorem{propn}{Proposition}
\newtheorem*{rmk*}{Remark}

\DeclareMathOperator{\dv}{\mbox{div}}
\DeclareMathOperator{\ord}{\mbox{ord}}
\begin{document}
\title[$q$-weights of branch points on superelliptic curves]{Higher-order Weierstrass weights of branch points on superelliptic curves}
\author[Caleb M. Shor]{Caleb McKinley Shor}
\address{Dept. of Mathematics, Western New England University, Springfield, MA 01119}
\email{cshor@wne.edu}
\begin{abstract}
In this paper we consider the problem of calculating the higher-order Weierstrass weight of the branch points of a superelliptic curve $C$.  For any $q>1$, we give an exact formula for the $q$-weight of an affine branch point. We also find a formula for the $q$-weight of a point at infinity in the case where $n$ and $d$ are relatively prime.  With these formulas, for any fixed $n$, we obtain an asymptotic formula for the ratio of the $q$-weight of the branch points, denoted $BW_q$, to the total $q$-weight of points on the curve:
\[ \liminf_{d\to\infty}\frac{BW_q}{g(g-1)^2(2q-1)^2}\geq \frac{n+1}{3(n-1)^2(2q-1)^2},\] with equality when the limit is taken such that $\gcd(n,d)=1$.
\end{abstract}

\keywords{Weierstrass points, higher-order Weierstrass points, superelliptic curves, branch points, numerical semigroups.}

\subjclass[2010]{Primary 14H55, 11G30}
\maketitle
\section{Introduction}
Let $q\in\mathbb{N}$.  A $q$-Weierstrass point (or higher-order Weierstrass point) is a point $P$ on a curve for which there exist holomorphic $q$-differentials that have higher than expected orders of vanishing at $P$.  Each $q$-Weierstrass point has an associated $q$-weight, denoted $w^{(q)}(P)$, which measures how much higher than expected those orders of vanishing are.  A curve of genus $g\geq2$ has finitely many $q$-Weierstrass points.

Importantly, the $q$-weight of a point is invariant under automorphism.  Thus, higher-order Weierstrass points are important in the study of automorphisms of algebraic curves.  For instance, Lewittes showed in \cite{Lewittes63} that if an automorphism has at least five fixed points, then all of its fixed points are $1$-Weierstrass points.  Further, Mumford, in \cite{Mumford99}, has suggested that $q$-Weierstrass points on an algebraic curve are analogous to $q$-torsion points on an elliptic curve.  For more on the history of Weierstrass points, we refer the reader to \cite{DelCentina2008}. For background material of Weierstrass points specifically on superelliptic curves, see \cite{ShaskaShor2015Advances}.

Let $C$ be a curve of genus $g\geq2$ of the form $y^n=f(x)$ where $f(x)$ is a separable polynomial of degree $d>n\geq2$.  Such a curve is said to be superelliptic.  In the cover $C\to\mathbb{P}^1$, the points above the roots of $f(x)$ are branch points.  If $n\nmid d$, the point (or points) above $\infty$ in the nonsingular model of $C$ is also a branch point.  One can show that each branch point is a $q$-Weierstrass point for all $q$; in the case where $n=2$, the branch points are exactly the $1$-Weierstrass points.  Let $B$ be an affine branch point on $C$ and, if $n\nmid d$, $P^\infty$ a nonsingular branch point at infinity. In \cite[Theorem 8]{Towse96}, Towse calculated the $1$-weight of the branch points (affine and at infinity) as a function of $n$ and $d$.  In the case of $\gcd(n,d)=1$, he found that \[w^{(1)}(B)=\dfrac{g(n+1)(d-7)}{12}+(d-1)\sum_{j=1}^{n-1}\left\{\dfrac{-dj}{n}\right\}j,\] where $\{x\}$ denotes the fractional part of $x$, and \[w^{(1)}(P^\infty)=\dfrac{(n^2-1)(d^2-1)}{24}-g.\]  Given that the total $1$-weight of points on a curve of genus $g$ is $g^3-g$, he was able to calculate the fraction that the branch points' $1$-weight (denoted $BW$) accounted for as \[ \lim_{d\to\infty} \dfrac{BW}{g^3-g}\geq\dfrac{n+1}{3(n-1)^2},\] with equality when the limit is taken over integers $d$ such that $\gcd(n,d)=1$.

The goal of this paper is to extend Towse's results to higher-order Weierstrass weights of branch points on a superelliptic curve.  To achieve this, we first produce a basis for the space of holomorphic $q$-differentials on a superelliptic curve $C$.  A common approach to calculate the $q$-weight is to work with the Wronskian of this basis, a method first described by Hurwitz in \cite{Hurwitz1893}. However, we take a different approach, instead using results from numerical semigroups and non-representable numbers. In this way we obtain a formula for the $q$-weight of an affine branch point as a function of $n$, $d$, and $q$.  The main result is Theorem~\ref{thm:main-result}.  In particular, when $\gcd(n,d)=1$ we find for $q\geq2$ \[w^{(q)}(B)=\dfrac{g(n+1)(d-7)}{12}+g+(d-1)\sum_{j=1}^{n-1}\left\{-\dfrac{(d+1)q+dj}{n}\right\}j.\]  As for the points at infinity, with two examples, we show that if $\gcd(n,d)>1$, one cannot get a formula for $w^{(q)}(P^\infty)$ based only on $n$, $d$, and $q$.  However, if $\gcd(n,d)=1$ and $q\geq2$, then in Theorem~\ref{thm:pt-at-infinity-weight} we have \[w^{(q)}(P^\infty)=\dfrac{(n^2-1)(d^2-1)}{24}.\]  Since the total $q$-weight of points on a curve of genus $g$ is $g(g-1)^2(2q-1)^2$ for $q\geq2$, we get an similar asymptotic result in Proposition~\ref{propn:proportion-of-branch-q-weight}  for the proportion of branch points' $q$-weight (denoted $BW_q$).  We find \[ \lim_{d\to\infty} \dfrac{BW_q}{g(g-1)^2(2q-1)^2} \geq\dfrac{n+1}{3(n-1)^2(2q-1)^2},\] with equality when the limit is taken over integers $d$ such that $\gcd(n,d)=1$.

The paper is organized as follows.  
In Section~\ref{sec:prelim}, we provide some background material on calculating the $q$-Weierstrass weight of a point.  We also include some notation and results for non-representable integers in numerical semigroups with two generators.  In Section~\ref{sec:basis}, we find a basis for the space of holomorphic $q$-differentials on a curve $C$ by presenting a set of linearly independent holomorphic differentials and then counting them to make sure there are as many as the Riemann-Roch theorem predicts.  In Section~\ref{sec:main-theorem}, we have our main results.  We find a formula for the $q$-weight of an affine branch point, and we use that to derive a few corollaries for specific cases of $n$ and $d$.  We also note that, for given $n$ and $d$, the $q$-weight of a branch point depends only on the value of $q$ modulo $n$. We also give examples to show that if $\gcd(n,d)>1$ the $q$-weight of a point at infinity cannot necessarily be determined just by knowing $n$ and $d$.  If $\gcd(n,d)=1$, however, we can calculate the weight.  In both cases, we obtain some asymptotic results about the proportion of $q$-weight that the branch points contain.  

\section{Preliminaries and notation}\label{sec:prelim}
\subsection{$q$-Weierstrass points}
In this paper, we will follow the approach given in \cite[Section 2]{ShaskaShor16}.  We describe the notation and major results for calculating weights of $q$-Weierstrass points here.

Let $k$ be an algebraically closed field,  $C$ be a non-singular projective curve over $k$ of genus $g\geq2$, and  $k(C)$ its function field.  For any $f\in k(C)$, let $\dv (f)$ denote the divisor associated to $f$.  For any divisor $D=\sum_P n_P P$ and any point $P$, let $\nu_P(D)=n_P$, and let $\ord_P(f)=\nu_P(\dv (f))$.

For any $q\in\mathbb{N}$, let $H^0(C,(\Omega^1)^q)$ be the $\mathbb{C}$-vector space of holomorphic $q$-differentials on $C$, a space of dimension 
\begin{equation*}\label{d_q}
d_q = 
\begin{cases} 
g           & \text{if } q=1, \\ 
(g-1)(2q-1) & \text{if } q\geq2.
\end{cases}
\end{equation*}

For $P$ a degree 1 point on $C$, consider a basis $\{\psi_1,\dots,\psi_{d_q}\}$ of $H^0(C,(\Omega^1)^q)$ where
\[ \ord_P(\psi_1)<\ord_P(\psi_2)<\cdots<\ord_P(\psi_{d_q}).\]
The $q$-weight of $P$ is \[w^{(q)}(P)=\sum_{i=1}^{d_q} \ord_P(\psi_i) - \sum_{j=0}^{d_q-1} j.\] We call the point $P$ a \textit{$q$-Weierstrass point} if $w^{(q)}(P)>0$.

\begin{propn}\cite[III.5.10]{FarkasKra92} \label{prop:total-q-weight}
Let $C$ be a curve of genus $g\geq2$ and let $q\geq1$.  Then the total $q$-weight of points on $C$ is \[\sum\limits_{P\in C} w^{(q)}(P) = (g-1)d_q(2q-1+d_q)= \begin{cases} g^3-g & \text{if $q=1$}, \\ g(g-1)^2(2q-1)^2 & \text{if $q\geq2$}.\end{cases}\]
\end{propn}

\subsection{Non-representable numbers}
For notation, let $\mathbb{N}_0$ be the set of non-negative integers.  Let $a,b\in\mathbb{N}$ and consider the set \[R(a,b)=\left\{ ax+by : x,y\in\mathbb{N}_0\right\}.\] Elements of $R(a,b)$ are called \textit{$(a,b)$-representable numbers}.  The complement of $R(a,b)$ in $\mathbb{N}_0$, denoted $NR(a,b),$ is the set of  \textit{$(a,b)$-representable numbers}.  When there is no confusion, we will omit the $(a,b)$ and simply refer to these numbers as representable or non-representable.

The problem of calculating the cardinality of $NR(a,b)$ dates to the late 19th century in \cite{Sylvester1882}.  Clearly, if $\gcd(a,b)>1$ then $NR(a,b)$ is an infinite set.  It is straightforward to show that the converse is true too.  For example, see \cite[Theorem 1.0.1]{RamirezAlfonsin05} for two proofs of the following result.
\begin{lem}
For $a,b\in\mathbb{N}$, if $\gcd(a,b)=1$ then $NR(a,b)$ is a finite set.
\end{lem}

For the rest of this section, we will assume $\gcd(a,b)=1$. Thus, $NR(a,b)$ is finite and so we can compute its cardinality.
\begin{propn}\cite[Page 134]{Sylvester1882}\label{propn:sylvester}
For $a,b\in\mathbb{N}$ with $\gcd(a,b)=1$, \[|NR(a,b)|=\frac{(a-1)(b-1)}{2}.\]\end{propn}

This result is important in the theory of algebraic curves.  Suppose a plane curve $C$ is given by the affine equation \[\alpha_{a,0}x^a + \alpha_{0,b}y^b  + \sum\limits_{i,j}\alpha_{i,j}x^iy^j = 0\] for constants $\alpha_{i,j}$ with $\alpha_{a,0}\cdot\alpha_{0,b}\neq0$ and where the summation is over non-negative $i,j$ such that $aj+bi<ab$.  Such a curve is called a $C_{a,b}$ curve.  These curves can be seen as a generalization of elliptic and hyperelliptic curves in Weierstrass form.  With the Riemann-Roch Theorem, if the affine part of the curve is non-singular then one can show that the genus of such a curve is exactly $(a-1)(b-1)/2$, the cardinality of $NR(a,b)$.  For details, see \cite{Miura93} or \cite{Shor11}.

For the purposes of this paper, we will also need to know the sum of the elements of $NR(a,b)$.  This problem was solved in \cite{BrownShiue93} using generating functions.

\begin{propn}\label{propn:sum-of-gaps}
For $a,b\in\mathbb{N}$ with $\gcd(a,b)=1$, \[\sum\limits_{n\in NR(a,b)} n = \dfrac{(a-1)(b-1)(2ab-a-b-1)}{12}.\]
\end{propn}
This result was generalized to a formula for the sum of the $m$th powers of elements of $NR(a,b)$.  See \cite{Rodseth93} or \cite{Tuenter06} for details.

\section{A basis of holomorphic $q$-differentials}\label{sec:basis}
In this section, we give a basis for the space of holomorphic $q$-differentials on a superelliptic curve $C$.  The main result of this section is as follows:

\begin{thm*}
Let $C$ be a curve of genus $g\geq2$ given in affine coordinates by $y^n=f(x)$, for $f(x)$ a separable polynomial of degree $d>n\geq2$.  For $q\geq1$, let $H^0(C,(\Omega^1)^q)$ be the space of holomorphic $q$-differentials on $C$. Let \[\mathfrak{B}_{n,d,q}=\left\{x^iy^j\left(\dfrac{\mathrm{d}x}{y^{n-1}}\right)^q : 0\leq i, \, 0\leq j<n, \, ni+dj\leq(2g-2)q\right\}.\]  Then $\mathfrak{B}_{n,d,q}$ is a basis for $H^0(C,(\Omega^1)^q)$.
\end{thm*}

In order to prove this, we need the following results (based on the work in \cite{Towse96}) and a useful lemma.

Let $G=\gcd(n,d)$.  For $g$ the genus of $C$, we have $2g-2=nd-n-d-G$.  Let $\{\alpha_1,\dots,\alpha_d\}$ denote the $d$ distinct roots of $f(x)$ and let $B_i=(\alpha_i, 0)$ for $i=1,\dots,d$.  For each non-root $\omega$ of $f(x)$, let $P_1^\omega,\dots,P_n^\omega$ denote the $n$ points on $C$ over $x=\omega$.  And let $P_1^\infty,\dots,P_G^\infty$ denote the $G$ points over $\infty$ in the non-singular model of $C$. 
One then has the following principal divisors.
\small
\begin{itemize}
\item $\dv (y)=\displaystyle \sum_{j=1}^{d} B_j  - \dfrac{d}{G}\sum_{m=1}^G P_m^\infty$,
\item $\dv (x- \alpha_i) = \displaystyle nB_i - \dfrac{n}{G}\sum_{m=1}^G P_m^\infty$,
\item $\dv (x-\omega) = \displaystyle \sum_{j=1}^n P_j^\omega - \dfrac{n}{G}\sum_{m=1}^G P_m^\infty$.
\item $\dv (\mathrm{d}x)=\displaystyle (n-1) \sum_{j=1}^{d} B_j  - \left(\dfrac{n}{G}+1\right)\sum_{m=1}^G P_m^\infty$,
\end{itemize}
\normalsize
\noindent From these, we see that $\dv ((\mathrm{d}x/y^{n-1})^q) = \frac{(2g-2)q}{G}\sum_{m=1}^G P_m^\infty$.

For integers $i,j$, let $f_{i,j}=x^iy^j(\mathrm{d}x/y^{n-1})^q$.  We want to find conditions on $i$ and $j$ such that $f_{i,j}$ is a holomorphic $q$-differential.  Note that $f_{i,j}$ can have poles only at the points above $\infty$ if $i,j\geq0$.  In that situation, we find
\[ \ord_{P_m^\infty}(f_{i,j}) = \dfrac{(2g-2)q-(ni+dj)}{G}\] for each pair $(i,j)$ and each $m$.  Hence $f_{i,j}$ is a holomorphic $q$-differential as long as $i\geq0, j\geq0,$ and $ni+dj\leq(2g-2)q$.

\begin{lem}\label{lem:exceptional-cases-n-d-q}
Let $n,d,q\in\mathbb{Z}$ with $2\leq n<d$ and $q\geq2$.  As above, let $G=\gcd(n,d)$ and $2g-2=nd-n-d-G$.
For all but finitely many triples $(n,d,q)$, one has \[(2g-2)q-d(n-1)\geq0.\]  The exceptional cases are $(n,d,q)\in\{(2,5,2), (2,6,2)\}.$
\end{lem}
\proof
First, note that $(2g-2)q-d(n-1)=(nd-n-d-G)q-d(n-1)=d(n-1)(q-1)-q(n+G)$.  Thus, to show our desired inequality, it is equivalent to show \[d\geq \left(\dfrac{q}{q-1}\right)\left(\dfrac{n+G}{n-1}\right).\]  For notation, let $h(q,n)=\left(\frac{q}{q-1}\right)\left(\frac{n+G}{n-1}\right).$ We aim to show $d\geq h(q,n)$.

For $q\geq2$, the maximum value of $q/(q-1)$, which occurs when $q=2$, is 2.  For $n\geq2$, the maximum value of $(n+G)/(n-1)$ occurs when $G$ is largest, so when $G=n$, and when $n=2$.  The maximum value is 4.  Thus $h(q,n)\leq2\cdot4=8$ so $d\geq h(q,n)$ for all $d\geq8$.

Now we consider cases of $n$.  If $n\geq4$, then $G=n$ or $G<n$.  If $G=n$, then $n|d$, so $d\geq2n\geq8$, so $d\geq h(q,n)$.  If $G<n$, then $G\leq n/2$, so $h(q,n)\leq 2(n+n/2)/(n-1)=3+3/(n-1)\leq4$.  Since $d>n$, we have $d>n\geq4\geq h(q,n)$, as desired.

If $n=3$, then $G=1$ or $G=3$.  If $G=1$, then $h(q,3)\leq 2\cdot2=4$.  In this case, since $d>n=3$, we have $d\geq4$, so $d\geq h(q,3)$.  If $G=3$, then $h(q,3)\leq 2\cdot3=6$.  In this case, $d\geq2n=6\geq h(q,3)$, as desired.

If $n=2$, then $G=1$ or $G=2$.  Note that to have $g\geq2$, we only consider $d\geq5$.  If $G=1$ and $q=2$, then $h(2,2)=2\cdot3=6$.  Since $G=1$, $d$ is odd, so $d\geq h(2,2)$ for all $d$ except for $d=5$.  If $G=1$ and $q\geq3$, then $h(q,2)\leq (3/2)\cdot3=9/2$, so $d\geq h(q,2)$ for all $d\geq5$.  If $G=2$ and $q=2$, then $h(2,2)=2\cdot4=8$.  Since $G=2$, $d$ is even, so $d\geq h(2,2)$ for all $d$ except for $d=6$.  If $G=2$ and $q\geq3$, then $h(q,2)\leq (3/2)\cdot4$, so $d\geq h(q,2)$ for all $d\geq6$.

Thus, the only exceptional cases are $(n,d,q)=(2,5,2)$ or $(2,6,2)$.  For all other triples, we find that $d\geq h(q,n)$, or, equivalently, that \[(2g-2)q-d(n-1)\geq0,\] as desired.
\qed

\begin{thm}
Let $C$ be a curve of genus $g\geq2$ given in affine coordinates by $y^n=f(x)$, for $f(x)$ a separable polynomial of degree $d>n$.  For $q\geq1$, let $H^0(C,(\Omega^1)^q)$ be the space of holomorphic $q$-differentials on $C$. Let \[\mathfrak{B}_{n,d,q}=\left\{x^iy^j\left(\dfrac{\mathrm{d}x}{y^{n-1}}\right)^q : 0\leq i, \, 0\leq j<n, \, ni+dj\leq(2g-2)q\right\}.\]  Then $\mathfrak{B}_{n,d,q}$ is a basis for $H^0(C,(\Omega^1)^q)$.\end{thm}

\proof The $q=1$ case, in a slightly different form, is proved in \cite{Towse96}.  For completeness, we will first prove the $q\geq2$ case here and then adapt our argument to cover the $q=1$ case.

Suppose $q\geq2$.  With the restriction that $0\leq j<n$, we see that these holomorphic $q$-differentials are linearly independent.  We therefore need to show that $\left\vert\mathfrak{B}_{n,d,q}\right\vert=d_q=(2q-1)(g-1)$.

We first consider the case where $(n,d,q)\not\in\{(2,5,2),(2,6,2)\}$ and let $\mathfrak{B}=\mathfrak{B}_{n,d,q}$. 
Note that we require $i\geq0$ and $ni+dj\leq(2g-2)q$, so \[0\leq i\leq \left\lfloor\dfrac{(2g-2)q-dj}{n}\right\rfloor.\]  For each $j=0,\dots,n-1$, we have $(2g-2)q-dj \geq (2g-2)q-d(n-1),$ which is non-negative by Lemma~\ref{lem:exceptional-cases-n-d-q}.  (This is why we handle the $(2,5,2)$ and $(2,6,2)$ cases separately.)  Thus, to calculate the number of pairs $(i,j)$ in $\mathfrak{B}$, we will let $j$ go from $0$ to $n-1$ and count the number of indices $i$ that correspond to each $j$ value.  I.e., \[\left\vert\mathfrak{B}\right\vert = \displaystyle\sum_{j=0}^{n-1} \left(1+ \left\lfloor\frac{(2g-2)q-dj}{n}\right\rfloor \right).\]
Since $\lfloor x\rfloor = x-\{x\}$, we simplify the sum to get \[\left\vert\mathfrak{B}\right\vert = n+(2g-2)q-\dfrac{d(n-1)}{2}-\sum_{j=0}^{n-1}\left\{\frac{(-d-G)q-dj}{n}\right\}.\]  

Now, we consider cases of $G$. If $G=n$, then $n|d$, so $n|((-d-G)q-dj)$, so each term in the summation is $0$.  Note that $n-d(n-1)/2=-(nd-d-2n)/2=-(nd-n-d-G)/2=-(g-1)$.  Then $\left\vert\mathfrak{B}\right\vert=(2g-2)q-(g-1)=(2q-1)(g-1)=d_q$, as desired.

Next, suppose $G\neq n$.  Let $n'=n/G$ and $d'=d/G$.  Dividing the numerator and denominator by $G$, the summation equals $\sum_{j=0}^{n-1}\left\{\frac{(-d'-1)q-d'j}{n'}\right\}.$  Since $\gcd(n',d')=1$, as $j$ goes from $0$ to $n'-1$ modulo $n'$, the numerators are distinct modulo $n'$ and therefore in every congruence class exactly once modulo $n'$.  Since $n/n'=G$, this summation equals $G \sum_{k=0}^{n'-1}\frac{k}{n'} = G(n'-1)/2$.  All together, \[\left\vert\mathfrak{B}\right\vert=n+(2g-2)q-d(n-1)/2-(n-G)/2.\]  I.e. $\left\vert\mathfrak{B}\right\vert=n+2q(g-1)-(1/2)(nd-d+n-G) = 2q(g-1)-(1/2)(nd-d-n-G) = (2q-1)(g-1)=d_q$, as desired.

To complete the proof for $q\geq2$, we consider the exceptional cases.  Suppose $(n,d,q)=(2,5,2)$, so $g=2$ and $d_2=3$.  Then
\[\mathfrak{B}_{2,5,2}=\left\{x^iy^j(\mathrm{d}x/y)^2 : i\geq0, 0\leq j<2, 2i+5j\leq 4.\right\}\]  Thus, $\mathfrak{B}_{2,5,2}=\{ (\mathrm{d}x/y)^2, \, x(\mathrm{d}x/y)^2, \, x^2(\mathrm{d}x/y)^2\}$, so $\left\vert\mathfrak{B}_{2,5,2}\right\vert=3=d_2$, as desired.

Suppose $(n,d,q)=(2,6,2)$, so $g=2$ and $d_2=3$.  Then
\[\mathfrak{B}_{2,6,2}=\left\{x^iy^j(\mathrm{d}x/y)^2 : i\geq0, 0\leq j<2, 2i+6j\leq 4.\right\}\]  Thus, $\mathfrak{B}_{2,6,2}=\{ (\mathrm{d}x/y)^2, \, x(\mathrm{d}x/y)^2, \, x^2(\mathrm{d}x/y)^2\}$, so $\left\vert\mathfrak{B}_{2,6,2}\right\vert=3=d_2$, as desired.

Now, suppose $q=1$.  Following the approach above, given $j$ we need integers $i$ such that \[0 \leq i \leq \left\lfloor\dfrac{(2g-2)-dj}{n}\right\rfloor.\]  If $j\leq n-2$, then $(2g-2)-dj\geq(2g-2)-d(n-2)=d-n-G\geq0$ since $d\geq n$.  If $j=n-1$, then $(2g-2)-d(n-1)=-n-G<0$, so there are no such $i$.  Thus, our summation for $\left\vert\mathfrak{B}\right\vert$ ends at $j=n-2$ instead of $j=n-1$.  Since we have a formula for the summation above, we can subtract the $j=n-1$ term out front to get \[\left\vert\mathfrak{B}\right\vert = -\left(1+\left\lfloor\dfrac{-n-G}{n}\right\rfloor\right)+ \sum\limits_{j=0}^{n-1}\left(1+ \left\lfloor\frac{(2g-2)-dj}{n}\right\rfloor \right).\]  Since $q=1$ the summation equals $g-1$, so $\left\vert\mathfrak{B}\right\vert=-(1-2)+(g-1)=g=d_1$, as desired.
\qed

\section{Weights of branch points}\label{sec:main-theorem}
In this section, we use the bases we found in the previous section to calculate the $q$-weight of the affine branch points and, in the case that $\gcd(n,d)=1$, the point at infinity.

\subsection{Weights of affine branch points}
Suppose $q\geq2$.  For $C$ given by $y^n=f(x)$ with $f(x)$ separable of degree $d$, let $\alpha$ be a root of $f(x)$.  Then $B=(\alpha,0)$ is an affine branch point of $C$.  Note that we can replace $x$ by $(x-\alpha)$ in our basis $\mathfrak{B}_{n,d,q}$ to produce a new basis $\mathfrak{B}_{n,d,q,\alpha}$.  That is, \[\mathfrak{B}_{n,d,q,\alpha}=\{(x-\alpha)^iy^j(\mathrm{d}x/y^{n-1})^q : i\geq0, \, 0\leq j<n, \, ni+dj\leq (2g-2)q\}\] is a basis for $H^0(C,(\Omega^1)^q)$.

Let $f_{i,j,\alpha}=(x-\alpha)^iy^j(\mathrm{d}x/y^{n-1})^q\in\mathfrak{B}_{n,d,q,\alpha}$.  Then \[\nu_{B}(f_{i,j,\alpha}) = ni+j. \]  Since $0\leq j<n$, these valuations are all different, and thus \[w^{(q)}(B)=\sum_{(i,j)\in S}(ni+j) - \sum_{k=0}^{d_q-1}k,\] where $S=\{(i,j)\in\mathbb{Z}^2 : i\geq0, \, 0\leq j<n, \, ni+dj\leq (2g-2)q\}$.  We rewrite this as $w^{(q)}(B)=W_1 - W_2 - W_3$ where 
\begin{equation}\label{eqn:W1-W2-W3}W_1=\sum_{(i,j)\in S}(ni+dj), \quad W_2=(d-1)\sum_{(i,j)\in S} j, \quad W_3=\sum_{k=0}^{d_q-1}k.\end{equation}  
We have \[W_3 = \frac{(d_q-1)(d_q)}{2}=\dfrac{1}{2}\left((2g-2)^2q^2 + (2g-2)(1-2g)q + g(g-1)\right).\]
We will evaluate $W_1$ and $W_2$ with the following propositions.

\begin{propn}\label{prop:W1}Let $n,d,q\in\mathbb{N}$ such that $n<d$ and $q\geq2$. 
Then \[W_1 = 2(g-1)^2q^2 + (g-1)Gq+\dfrac{G^2-1-(n-1)(d-1)(2nd-n-d-1)}{12}.\]
\end{propn}

We will first sketch the proof in the situation where $\gcd(n,d)=1$.  Afterward, we will prove the theorem for any gcd.

When $\gcd(n,d)=G=1$, for $(i,j)\in S$, the terms $ni+dj$ are distinct integers from $0$ to $(2g-2)q$.  From Proposition \ref{propn:sylvester}, since $(2g-2)q\geq nd-n-d$ (by Lemma~\ref{lem:mk-bound} below), all of the $(n-1)(d-1)/2$ $(n,d)$-non-representable integers are in that interval.  The sum of the non-representable integers, as is given in Proposition~\ref{propn:sum-of-gaps}, is $(n-1)(d-1)(2nd-n-d-1)/12$.  Thus, if $\gcd(n,d)=1$, we add up all of the integers from $0$ to $(2g-2)q$ and subtract off the non-representable integers to get $W_1 = (2g-2)q((2g-2)q+1)/2 - (n-1)(d-1)(2nd-n-d-1)/12$.

If $\gcd(n,d)>1$, then the terms $ni+dj$ are no longer distinct, so we need to evaluate the sum more carefully.

\proof Let $G=\gcd(n,d)$.  First, we observe that $W_1=G\sum_{(i,j)\in S} (n'i+d'j)$ for $n'=n/G$ and $d'=d/G$. Note that $\gcd(n',d')=1$.  For $k$ from $0$ to $G-1$, let \[ S_k = \{(i,j) : i\geq0, \, kn'\leq j< (k+1)n', \, ni+dj\leq (2g-2)q\}.\] In particular, $S$ is the disjoint union of the sets $S_k$.  Let $W_{1,k}=\sum_{(i,j)\in S_k}(n'i+d'j).$  Then \[W_1 = G\sum\limits_{k=0}^{G-1} W_{1,k}.\]

Letting $j'=j-kn'$, we rewrite $S_k$ as \[ S_k = \{(i,j'+kn') : i\geq0, \, 0\leq j'<n', \, ni+dj'\leq(2g-2)q-n'dk\}\] and dividing the last inequality through by $G$ we obtain \[ S_k = \{(i,j'+kn') : i\geq0, \, 0\leq j'<n', \, n'i+d'j'\leq\dfrac{(2g-2)}{G}q-n'd'k\}.\]  Let $m_k=\frac{(2g-2)}{G}q-n'd'k$, the upper bound in $S_k$.  The following lemma will allow us to conclude that all of the $(n', d')$-non-representable integers are less than $m_k$. 

\begin{lem}\label{lem:mk-bound}
Let $m_{k}=\frac{(2g-2)}{G}q-n'd'k.$  Then $m_{k}\geq n'd'-n'-d'$ for all $n,d,q\in\mathbb{N}$ with $0\leq k\leq G-1$, $n<d$, $g\geq2$, and $q\geq2$.
\end{lem}
\proof
First, note that for $0\leq k\leq G-1$, $m_k=\frac{(2g-2)}{G}q-n'd'k\geq\frac{(2g-2)}{G}q-n'd'(G-1)$.  So we need to show $\frac{(2g-2)}{G}q-n'd'(G-1)\geq n'd'-n'-d'$, which is equivalent to showing $(2g-2)q\geq nd-n-d$.  Since $nd-n>nd-n-d$, by Lemma~\ref{lem:exceptional-cases-n-d-q}, we have $(2g-2)q\geq nd-d>nd-n-d$ for all $(n,d,q)$ combinations except $(2,5,2)$ and $(2,6,2)$.

We compute the exceptional cases separately.  If $(n,d,q)=(2,5,2)$, then $g=2$ and $(2g-2)q=4\geq3=nd-n-d$.  If $(n,d,q)=(2,6,2)$, then $g=2$ and $(2g-2)q=4\geq4=nd-n-d$.  Thus, the bound holds for the exceptional cases as well.
\qed

For $S_k$, since we are considering $i\geq0$ and $0\leq j'<n'$, and since $m_k\geq n'd'-n'-d'$ for all $k$, our ordered pairs $(i, j'+kn')$ are in one-to-one correspondence with the $(n',d')$-representable numbers in the interval $[0, m_k]$.  And since $m_k\geq n'd'-n'-d'$, all of the $(n'-1)(d'-1)/2$ $(n',d')$-non-representable numbers are in this interval as well. Thus $S_k$ contains $\left\vert S_k\right\vert=m_k+1-(n'-1)(d'-1)/2$ ordered pairs.

Then \begin{align*}
W_{1,k} &= \sum\limits_{(i,j)\in S_k} (n'i+d'j) \\
&= \sum\limits_{(i,j'+kn')\in S_k} \left(n'i+d'j'+n'd'k\right) \\
&= n'd'k\cdot\left\vert S_k\right\vert+\sum\limits_{(i,j'+kn')\in S_k} \left(n'i+d'j'\right).\\
\end{align*}
The summation is the sum of the $(n',d')$-representable numbers from $0$ to $m_k$.  We calculate this by summing all of the integers from $0$ to $m_k$ and subtracting the $(n',d')$-non-representable integers, which all lie in this interval. Using Proposition~\ref{propn:sum-of-gaps}, the summation is $m_k(m_k+1)/2 - (n'-1)(d'-1)(2n'd'-n'-d'-1)/12$.  Thus,
\begin{align*}
W_{1,k} & = n'd'k\left(m_k+1-(n'-1)(d'-1)/2\right)& \\ 
& + m_k(m_k+1)/2-(n'-1)(d'-1)(2n'd'-n'-d'-1)/12, &
\end{align*}
so
\begin{align*}
W_1 = &G\sum\limits_{k=0}^{G-1} \left[n'd'k\left(m_k+1-\dfrac{(n'-1)(d'-1)}{2}\right)+\dfrac{m_k(m_k+1)}{2}\right. & \\ 
& \hspace{30px} -\left.\dfrac{(n'-1)(d'-1)(2n'd'-n'-d'-1)}{12}\right].&
\end{align*}

To evaluate this sum, we need the following calculations which are straightforward to compute.

\begin{itemize}
\item $\sum\limits_{k=0}^{G-1} m_k = (2g-2)q-\frac{n'd'G(G-1)}{2}$.
\item $\sum\limits_{k=0}^{G-1} m_k^2 = \frac{(2g-2)^2}{G}q^2 - (2g-2)(G-1)d'n'q+\frac{d'^2n'^2(G-1)G(2G-1)}{6}$
\item $\sum\limits_{k=0}^{G-1} k m_k = (g-1)(G-1)q - \frac{d'n'(G-1)G(2G-1)}{6}$.
\end{itemize}

Simplifying the resulting expression, we find \[W_1 = \dfrac{(2g-2)^2}{2}q^2 + (g-1)Gq+\dfrac{G^2-1-(n-1)(d-1)(2nd-n-d-1)}{12},\] which completes the proof of Proposition~\ref{prop:W1}.
\qed

\begin{propn}\label{propn:W2}
Let $n,d,q\in\mathbb{N}$ such that $n<d$ and $q\geq2$.  Let \[D(a,b,c)=\sum_{j=0}^{c-1}\left\{\dfrac{a+bj}{c}\right\}j.\] 
Then \[W_2 = (d-1)\left((n-1)\left((g-1)q+\dfrac{-2nd+3n+d}{6}\right) - D(-(d+G)q,-d,n)\right).\]
\end{propn}

\proof
We will use Lemma~\ref{lem:exceptional-cases-n-d-q}, so we first assume $(n,d,q)\not\in\{(2,5,2),(2,6,2)\}$. For $W_2=(d-1)\sum_{(i,j)\in S} j$, we have \[W_2=(d-1)\sum_{j=0}^{n-1}\sum_{i=0}^{I_j}j,\] for $I_j=\left\lfloor\frac{(2g-2)q-dj}{n}\right\rfloor.$  By Lemma~\ref{lem:exceptional-cases-n-d-q}, $I_j\geq0$ so $W_2=(d-1)\sum_{j=0}^{n-1}(I_j+1)j$.  Since $\lfloor x\rfloor=x-\{x\}$, \[W_2=(d-1)\sum_{j=0}^{n-1}\left(\dfrac{(2g-2)q-dj}{n}-\left\{\dfrac{(2g-2)q-dj}{n}\right\}+1\right)j.\]  Note that $\left\{\frac{(2g-2)q-dj}{n}\right\}=\left\{\frac{(nd-n-d-G)q-dj}{n}\right\}=\left\{\frac{(-d-G)q-dj}{n}\right\}.$

Expanding out, we get
\begin{align*}
W_2 = & (g-1)(d-1)(n-1)q + \dfrac{(d-1)(n-1)}{6}(3n-d(2n-1)) & \\
& \hspace{50px} - (d-1)\sum_{j=0}^{n-1}\left\{-\dfrac{(d+G)q+dj}{n}\right\}j,& \end{align*} which can be rearranged to give the desired result.

Finally, if $(n,d,q)\in\{(2,5,2),(2,6,2)\}$, then $S=\{(0,0),(1,0),(2,0)\}$, and so $W_2=\sum_{(i,j)\in S} j = 0$.  We get the same value if we plug each these $(n,d,q)$ triples into the above formula for $W_2$.
\qed

\begin{rmk*}
The summation $D(a,b,c)$ is related to a Dedekind sum.  There is no closed form for such sums, though there is a reciprocity law.  For a general reference, see \cite{Rademacher72}. 
\end{rmk*}

Finally, we can combine and simplify $W_1-W_2-W_3$.  Note that the $q^2$ and $q$ terms (other than in the summation) cancel.  With further manipulation, we have our main result.

\begin{thm}\label{thm:main-result} Let $C$ be given in affine coordinates by $y^n=f(x)$ for $f(x)$ a separable polynomial of degree $d>n$.  Let $G=\gcd(n,d)$, and let $q\in\mathbb{Z}$ with $q\geq2$.  For any root $\alpha$ of $f(x)$, let $B=(\alpha,0)$ be a branch point.

The $q$-weight of $B$ is $w^{(q)}(B) = $
\begin{align*}
w^{(q)}(B) = & \dfrac{1}{24}\left((n-1)(d-1)(n+1)(d-7)+12g(G+1)+5(G^2-1)\right) & \\
& \hspace{50px} + (d-1)\cdot D(-(d+G)q,-d,n)&
\end{align*}
\end{thm}

Note that, for given values of $n$ and $d$, the $q$-weight of $B$ depends only on the value of $q$ modulo $n$.

We will give results for some combinations of $n$ and $d$ in the corollaries below.  First, we consider the case where $\gcd(n,d)=1$.

\begin{cor}If $\gcd(n,d)=1$, \[w^{(q)}(B) = \frac{g}{12}(n+1)(d-7)+g+(d-1)\cdot D(-(d+1)q,-d,n).\]
\end{cor}

Fix $n$ and $d$ (with any gcd).  If one varies $q$, then one sees the value of $w^{(q)}(B)$ depends only on the congruence class of $q$ modulo $n/G$.  Further, if $d\equiv-G\text{ (mod }n)$, then the summation term simplifies to $\sum_{j=0}^{n-1}\left\{\frac{Gj}{n}\right\}j,$ for which there is a closed form.

\begin{cor}
If $d\equiv-G\text{ (mod }n)$, then $w^{(q)}(B)$ doesn't depend on $q$.  In particular, 
\begin{align*}
w^{(q)}(B)=&\dfrac{1}{24}\biggl((n-1)(d-1)(n+1)(d-7)+12g(G+1)+5(G^2-1) & \\
&\hspace{50px}+ 2(d-1)(n-G)(3n+n'-2)\biggr).&\end{align*}
\end{cor}
\proof
The summation term is $\sum_{j=0}^{n-1}\left\{\frac{Gj}{n}\right\}j,=\sum_{j=0}^{n-1}\left\{\frac{j}{n'}\right\}j$.  Each $j$ can be written uniquely as $j=j'+kn'$ for $0\leq k<G$ and $0\leq j'<n'$.  Thus, the summation is $\sum_{k=0}^{G-1} \sum_{j'=0}^{n'-1} \left\{\frac{j'}{n'}\right\}j = \sum_{k=0}^{G-1} \sum_{j'=0}^{n'-1} \left(\frac{j'^2}{n'}+j'k\right)$, which simplifies to $(n-G)(3n+n'-2)/12$.
\qed

Combining the two corollaries above, we obtain the following.

\begin{cor}\label{cor:-1modn}
If $d\equiv-1\text{ (mod }n)$, then \[w^{(q)}(B) = \dfrac{g(n+1)(d+1)}{12} = \dfrac{(n^2-1)(d^2-1)}{24}\] for all $q\geq2$.
\end{cor}

\begin{cor}
If $n\mid d$, then \[w^{(q)}(B)=\frac{(n^2-1)(d^2-2d)}{24}.\]
\end{cor}
\proof
If $n\mid d$, then $G=n$ and $n\mid ((d+G)q+dj)$ for all $j$, so the summation is zero.  Since $2g-2=nd-n-d-n$, we have $g=\frac{(d-2)(n-1)}{2}$.  Plugging in, the result follows.
\qed

\subsection{Weights of points at infinity}
If $n\mid d$, then there are $n$ points at infinity in the smooth model of $C$, so these points are not branch points.  However, we can still investigate their $q$-weights.  If $\gcd(n,d)>1$, then we need to know more about $f(x)$ to determine $w^{(q)}(P_m^\infty)$.  We give a few examples to illustrate this.

In \cite{FarahatSakai11}, the authors consider curves of the form $y^2=f(x)=x^6+ax^4+bx^2+1$, where $a,b$ are parameters and $f(x)$ is separable.  In the non-singular models of these curves, there are $G=\gcd(n,d)=2$ points at infinity $P_1^\infty$ and $P_2^\infty$.  If $4b=a^2$, then $w^{(3)}(P_1^\infty)=w^{(3)}(P_2^\infty)=2$.  If $4b\neq a^2$, then $w^{(3)}(P_1^\infty)=w^{(3)}(P_2^\infty)=0$.  

In \cite[Lemma 4 and Proposition 3]{ShaskaShor16}, the authors consider hyperelliptic curves of genus 3 of the form $y^2=f(x)$ where $\deg(f)=8$.  In the non-singular models of these curves, there are $G=\gcd(n,d)=2$ points at infinity $P_1^\infty$ and $P_2^\infty$.  If $C$ is given by $y^2=x^8+x^6+16x^4+x^2+1$, then $w^{(2)}(P_1^\infty)=w^{(2)}(P_2^\infty)=1$. If $C$ is given by $y^2=x^8+x^4+1$, then $w^{(2)}(P_1^\infty)=w^{(2)}(P_2^\infty)=3$.

Thus, simply knowing $n$ and $d$ is not enough to calculate the $q$-weight of the points at infinity.  However, there are some cases where we can get a result.

First, if $d=n+1$, then the lone point at infinity is a nonsingular branch point, so it will have the same $q$-weight as the affine branch points.  By Corollary~\ref{cor:-1modn}, since $d\equiv-1\text{ (mod $n$)}$, $w^{(q)}(B)=\frac{(n^2-1)(d^2-1)}{24}$ for $q\geq2$, so we will have $w^{(q)}(P^\infty)=\frac{(n^2-1)(d^2-1)}{24}$ for $q\geq2$ as well.  This is a special case of the more general result when $\gcd(n,d)=1$.

\begin{thm}\label{thm:pt-at-infinity-weight}
Suppose $C$ is a curve of genus $g\geq2$ given by the affine equation $y^n=f(x)$ for $f(x)$ a separable polynomial of degree $d$ where $n<d$ and $\gcd(n,d)=1$.  Let $P_1^\infty$ be the lone point at infinity in the non-singular model of $C$.  Then \[w^{(q)}(P_1^\infty) = \begin{cases}\dfrac{g(n+1)(d+1)}{12}-g=\dfrac{(n^2-1)(d^2-1)}{24}-g & \text{if $q=1$,} \\ \dfrac{g(n+1)(d+1)}{12}=\dfrac{(n^2-1)(d^2-1)}{24} & \text{if $q\geq2$}.\end{cases}\]
\end{thm}
\proof
For $q=1$, the formula is given at the end of the proof of \cite[Theorem 8]{Towse96}.

For $q\geq2$ and $G=1$, let $\mathfrak{B}_{n,d,q}$ be as in Section~\ref{sec:basis}, and again let $S=\{(i,j)\in\mathbb{Z}^2 : i\geq0, 0\leq j<n, ni+dj\leq(2g-2)q\}$.  Then $f_{i,j}\in\mathfrak{B}_{n,d,q}$ if and only if $(i,j)\in S$.  Recall that $\ord_{P_m^\infty}(f_{i,j})=(2g-2)q-(ni+dj)$. These orders of vanishing are unique, so
\[w^{(q)}(P_m^\infty) = \left( \sum\limits_{(i,j)\in S} \ord_{P_m^\infty}(f_{i,j}) \right) - \sum\limits_{k=0}^{d_q-1} k.\]
Since $\left\vert S\right\vert=d_q$, \[w^{(q)}(P_m^\infty)=d_q(2g-2)q-\left( \sum\limits_{(i,j)\in S}(ni+dj) \right) - \dfrac{(d_q-1)d_q}{2}.\]  The summation, which we called $W_1$ in Equation~\ref{eqn:W1-W2-W3}, is evaluated in Proposition~\ref{prop:W1}.  Plugging this and $d_q$ in, the expression simplifies to $w^{(q)}(P_m^\infty)=\frac{(n^2-1)(d^2-1)}{24}$.
\qed

\subsection{Branch weight}
In the case where $\gcd(n,d)=1$, we can calculate the total $q$-weight of the branch points (both affine and at infinity) for $q\geq2$, which we denote $BW_q$.

\begin{cor}
Suppose $\gcd(n,d)=1$, so $g=\frac{(n-1)(d-1)}{2}$.  Then the total branch $q$-weight is given by 
$BW_q = d \cdot w^{(q)}(B) + w^{(q)}(P_1^\infty) =$
\[d\left( \frac{g}{12}(n+1)(d-7)+g+(d-1)\cdot D(-(d+1)q,-d,n)
\right) + g\frac{(n+1)(d+1)}{12}.\]  
Rewritten in terms of $g$, we get \[BW_q = \frac{n+1}{3(n-1)^2}\left(g^3-2g^2(n-1)-g(n-1)^2\right)+d(d-1)\cdot D(-(d+1)q,-d,n).\]\end{cor}

From Proposition~\ref{prop:total-q-weight}, we know the total weight of the $q$-Weierstrass points, for $q\geq2$, is $g(g-1)^2(2q-1)^2$.  We can now calculate the proportion of $q$-weight of the branch points.

\begin{propn}\label{propn:proportion-of-branch-q-weight}
Fix $n$ and let $q\geq2$.  Then \[\liminf\limits_{d\to\infty} \frac{BW_q}{g(g-1)^2(2q-1)^2} \geq \frac{n+1}{3(n-1)^2(2q-1)^2}.\]
If we restrict to values of $d$ that are relatively prime to $n$ then 
\[\lim\limits_{d\to\infty, (n,d)=1} \frac{BW_q}{g(g-1)^2(2q-1)^2} = \frac{n+1}{3(n-1)^2(2q-1)^2}.\]
\end{propn}
\proof
For general $n$ and $d$, since we do not have an exact formula for the $q$-weight of the points at infinity, we can only say $BW_q\geq d \cdot w^{(q)}(B)$.  Using the result from Theorem~\ref{thm:main-result}, since \[\sum_{j=0}^{n-1}\left\{-\frac{(d+G)q+dj}{n}\right\}j \leq \sum_{j=0}^{n-1} \frac{n-1}{n} j = \frac{(n-1)^2}{2},\] in terms of $d$, the dominant term of $d\cdot w^{(q)}(B)$ is $d^3\frac{(n-1)(n+1)}{24}$.  Since $g$ is on the order of $d(n-1)/2$, the dominant term of the denominator is $d^3\frac{(n-1)^3(2q-1)^2}{8}.$  The result follows.

For $\gcd(n,d)=1$, the lone point at infinity has weight $\frac{(d^2-1)(n^2-1)}{24}$.  Thus, the dominant term of $BW_q$ is precisely $d^3\frac{(n-1)(n+1)}{24}$, and we thus have an equality if we take a limit involving integers $d$ such that $\gcd(n,d)=1$.
\qed

\bibliography{main-bib-file}{}
\bibliographystyle{plain}
\end{document}